\documentclass[a4paper,11pt]{amsart}

\usepackage{amssymb,amsbsy,amsmath,amsfonts,amssymb,amscd}
\usepackage{latexsym}
\usepackage{graphics}
\usepackage{color}
\input xy
\xyoption{all}

\theoremstyle{plain}
\newtheorem{thm}{Theorem}[section]

\theoremstyle{definition}

\newtheorem{example}[thm]{Example}

\numberwithin{equation}{section}
\newcommand{\sR}{{\mathcal R}}

\newcommand{\CC}{\ensuremath{\mathbb{C}}}
\newcommand{\RR}{\ensuremath{\mathbb{R}}}
\newcommand{\ZZ}{\ensuremath{\mathbb{Z}}}

\newcommand\om{\omega}

\newcommand\ga{\gamma}

\def\eea{\end{eqnarray*}}
\def\bea{\begin{eqnarray*}}

\newcommand\dual{\mathrel{\raise3pt\hbox{$\underline{\mathrm{\thinspace d
\thinspace}}$}}}
\newcommand\qe{\ifhmode\unskip\nobreak\fi\quad $\Box$}       % box for QED

\def\BOX{\hfill\lower.5\baselineskip\hbox{$\Box$}}
\newtheorem{theo}{Theorem}[section]

\newtheorem{remarkk}[theo]{Remark}
\newenvironment{rem}{\begin{remarkk}\rm}{\end{remarkk}}

\newtheorem{prop}[theo] {Proposition}
\newtheorem{cor}[theo]{Corollary}

\newcommand{\Proof}{{\it Proof. }}

\title [$D_4$-HYperelliptic 3fold]{Hyperelliptic threefolds with group $D_4$, the dihedral group of order 8}

\author{Fabrizio Catanese and Andreas Demleitner}
\address {Lehrstuhl Mathematik VIII\\
Mathematisches Institut der Universit\"at Bayreuth\\
NW II,  Universit\"atsstr. 30\\
95447 Bayreuth}
\email{Fabrizio.Catanese@uni-bayreuth.de \newline \hspace*{2.2cm} Andreas.Demleitner@uni-bayreuth.de}

\thanks{
The present work took place in the realm of the 
 ERC Advanced grant n. 340258, `TADMICAMT' }

\date{\today}

\begin{document}

\maketitle

\begin{abstract}
We give a simple construction for the hyperelliptic threefolds with group $D_4$.

\end{abstract}

%\tableofcontents

\section*{Introduction}

A Generalized Hyperelliptic Manifold is the quotient $ X = T / G$ of a complex torus $T$ by the free action of a finite group $G$ which contains  no translations. 
We say that we have a Generalized 
Hyperelliptic Variety if moreover the torus $T$ is projective, i.e., it is an Abelian variety $A$.

 Recently D. Kotschick observed that the classification of Generalized Hyperelliptic Manifolds of complex dimension three  was not complete,
since the case where $G$ is the dihedral group $D_4$ of order $8$ was excluded (by H. Lange in  \cite{Lange}) but it does indeed occur.  Indeed 
 F.E.A. Johnson in the preprint \cite{Johnson} showed that 
 a construction due to Dekimpe, Ha\l enda and Szczepa\'nski of a flat manifold $M$  of real dimension 6 with holonomy equal to $D_4$ (see \cite{Dekimpe})
 would give the desired Manifold (which is projective, as remarked by Kotschick \footnote{personal communication to the second author}, being K\"ahler with second Betti number $=2$).
 
We describe  all such  examples explicitly, following the method of Lange, which was based on the classification of automorphisms of complex tori of dimension $2$ given by Fujiki  in \cite{Fujiki}.

The family we give is exactly the one obtained by taking all possible complex structures on the flat manifold $M$, and the upshot is that
all these hyperelliptic complex manifolds $X$ are quotients of the product of three elliptic curves by a translation of order 2. 

\section{The example}

Let $E, E'$ be any two elliptic curves, $$ E = \CC / ( \ZZ + \ZZ \tau ), \ \  E' = \CC / ( \ZZ + \ZZ \tau' ).$$

Set $$A' : = E \times E \times E' ,   A  : = A' /  \langle \om\rangle ,  \text{ where } \om : = (1/2, 1/2, 0).$$

\begin{thm}\label{example}

The Abelian variety $A$ admits a  fixed point free action of the dihedral group $$D_4 : = \langle r,s | r^4=1, s^2=1, (rs)^2 = 1\rangle,$$
such that $D_4$ contains no translations.

\end{thm} 

\Proof
Set, for $ z : = (z_1, z_2, z_3 ) \in A'$:

$$ r (z_1, z_2,z_3) : = (z_2, - z_1, z_3 + 1/4) = R (z_1, z_2,z_3) + (0,0,1/4)$$
$$ s (z_1, z_2, z_3) : = (z_2 + b_1 , z_1 + b_2, - z_3 ) = S (z_1, z_2,z_3) + (b_1, b_2, 0 ),  $$
$$\text{where }  b_1 := 1/2 + \tau /2 , \ b_2 := \tau/2.$$

{\bf Step 1.}
It is easy to verify that $r,R$ have order exactly $4$ on $A'$, and that $R (\om) =\om$, so that $r$ descends to an automorphism of $A$, of order exactly $4$.
Moreover, any power $r^j, \ 0 < j \leq 3$ acts freely on $A$, since the third coordinate of $r^j(z)$ equals $z_3 + j/4$. 

{\bf Step 2.}
$s^2(z) = z + \om$, since $b_1 + b_2 = 1/2$; moreover $S (\om) =\om$, hence $s$ descends to an automorphism of $A$ of order exactly $2$.

{\bf Step 3.} We have
$$ rs (z) = (z_1 + b_2, - z_2  - b_1, - z_3 + 1/4), $$ hence $$(rs)^2 (z) = ( z_1 + 2 b_2 , z_2, z_3) = z,$$
and we have an action of $D_4$ on $A$, since the respective orders of $r, s, rs$ are precisely $4,2,2$.

{\bf Step 4.}
We claim that also the symmetries in $D_4$ act freely on $A$ and are not translations. Since there are exactly two conjugacy classes of symmetries, those of $s$ and $rs$, it suffices to 
observe that  these two transformations are not translations: in the next step we show that they both act freely.

{\bf Step 5.}
It is rather immediate to see that $rs$ acts freely, since $rs (z) = z$ is equivalent to 
$$ (b_2, - 2 z_2 -b_1, - 2 z_3 + 1/4) $$
being a multiple of $\om$ in $ A'$. But this is absurd, since $2 \om = 0$, and $b_2 = \tau /2 \neq 0, 1/2$. 

The condition that $s$ acts freely, since $s(z) = z$ 
is equivalent to 
$$ \ga : = (z_2 - z_1 + b_1, z_1 -  z_2  + b_2, - 2 z_3 ) $$
being a multiple of $\om$ in $A'$. 

However the sum of the first two coordinates of multiples of $\om$ is $0$, and this equation is not satisfied by $\ga$, since $b_1 + b_2 = 1/2 \neq 0$.

\qed
\begin{cor}
The above family of hyperelliptic threefolds $X$ with group $D_4$  forms a  complete two dimensional family.
The K\"ahler manifolds with the same fundamental group as $X$ yield an open subspace of the Teichm\"uller space of $X$ 
parametrized by the  two halfplanes containing $\tau$, $\tau '$ respectively.

\end{cor}
\Proof
Take the above family of varieties $ X = A / G$ , where $G  = D_4$, and observe that, setting $U := \CC^3$,  $ H^1 (\Theta_A)^G = (U \otimes \overline{U}^{\vee})^G$
can be calculated as follows. We have $U = U_1 \oplus W$, where $U_1, W$ are real (self-conjugate) representations, $U_1$ is irreducible and $W$ is a character of $ G / \langle r\rangle$.

Hence $(U \otimes \overline{U}^{\vee})^G  = (U \otimes U^{\vee})^G = End (U_1 \oplus W)^G $ has dimension 2 by Schur's lemma.

Following Theorem 1 of \cite{Catanese-Corvaja}, and since, as we show below in proposition \ref{HT}, there is only one possible Hodge-Type, we conclude that the open subspace of the Teichm\"uller space of $X$ 
corresponding to K\"ahler manifolds is irreducible and equal to  the product of two halfplanes.

\qed

\begin{rem}
Let $$ \sR : = \sR_4 : = \ZZ [x] / (x^2 + 1)$$

be the 4th cyclotomic ring, also called the ring of Gaussian integers.

We denote by $\sigma$ the Galois involution sending $z = a  + x b \mapsto \sigma (z) := a - x b$.

We define, according to Dekimpe et al. (\cite{Dekimpe})  the following $\sR$-module:

$$ L : = \sR \oplus \sR \oplus (\ZZ e_5 \oplus  \ZZ e_6) = : L_1  \oplus L_2 \oplus L_3, $$ 
where the module $L_3$ is the trivial $\sR$-module.

The real torus $ T := (L \otimes \RR ) / L$ admits a free action of the dihedral group $D_4$,
defined as follows:

$$ r (z_1,z_2,z_3) = ( x z_1 + x/2, x z_2 + 1/2, z_3 - 1/4 e_5),$$
$$  s_2 (z_1,z_2,z_3) = ( x \sigma (z_1) , \sigma(z_2) + 1/2, - z_3).$$ 

It is easy to see that   the flat manifolds $A/D_4$ are the same as the flat manifolds  $T/ D_4$.

In fact, we define $ \om_1 , \om_2 \in \CC^3 $ as the vectors
 $$\om_1 : = (1/2, 1/2, 0), \om_2 := (\tau, 0,0).$$
 
  Then we set $L_1$ to be  the free $\sR$-module generated by $\om_1$, $L_2$   the free $\sR$-module generated by $\om_2$, 
  and $L_3$ the trivial $\sR$-module generated by $(0,0,1)$ and $(0,0,\tau')$.
  
  We have then  that $A$ is the quotient $\CC^3 / L$.

\end{rem}
\begin{prop}\label{HT}
The above $\sR$-module $L$ admits a unique Hodge-Type.

More precisely, the $D_4$-invariant complex structures form a $2$-dimensional complex family,
 obtained choosing  respective $1$ dimensional subspaces $U(1) \subset V(1) = V_3$
 and $U(i) \subset V(i)$, such that
 $$ (**)  \  \ U(1) \oplus \overline{U(1) } =  V(1), \  \  U(i) \oplus \overline{s U(i) } =  V(i) ,$$ 
 and defining $U(-i) : = S U(i)$.

\end{prop}
\Proof

Consider the complex vector space 
$$ V := L \otimes \CC = (L_1 \otimes \CC )\oplus (L_2\otimes \CC)  \oplus (L_3 \otimes \CC), $$
where we observe that each summand is stable by the action of $D_4$.

By looking at the eigenspace decomposition of $V$ with respect to the linear action of $r$, given by the diagonal matrix with entries $(x,x,1)$
we can decompose:
$$ V = (V_1(i) \oplus V_1(-i) ) \oplus (V_2(i) \oplus V_2(-i)) ) \oplus V_3.$$

The second summand is the conjugate of the first, the fourth is the conjugate of the third.
To get a free action of $D_4$ one must  give a Hodge decomposition
$$ V = U \oplus \overline{U}, $$
where the holomorphic subspace $U$ must be invariant by $R$, and must split as the sum of three eigenspaces for $R$:
$$  U = U(i) \oplus U(-i) \oplus U(1).$$

Let us look at $V_1(i)$: it is spanned by $(1-ix,0,0)$ and the linear part of $s$ sends it to $-i (1+ix)$, an element in $V_1(-i) $.

Similarly $S$ sends $V_2(i)$ to $V_2(-i)$, hence $$S V(i) = V(-i) = \overline {V(i)}.$$

Since $S$ preserves the complex structure $U$, we obtain that 
$$S (U(i)) = U(-i),$$ so all eigenspaces have dimension one.

We have already seen  that it must hold $U(-i)  = S U(i)$, from which follows that $S(U(i)) = U(-i)$.

The condition $V = U \oplus \overline{U}$ amounts then to the two properties $(**)$.
One can directly verify that the second  holds true  on some open set of the Grassmannian.

However, this also follows from the explicit description of our examples.

\qed

\end{document}